\documentclass[12pt,a4paper]{amsart}
\usepackage{amsfonts,amssymb}

\newcommand{\B}{\Bbb{B}}\newcommand{\C}{\Bbb{C}}\newcommand{\K}{\Bbb{K}}\newcommand{\bK}{\bar{\K}}
\newcommand{\mN}{\Bbb{N}}
\newcommand{\R}{\Bbb{R}}\newcommand{\N}{\Bbb{N}}
\newcommand{\li}{~\\ $\bullet$ }\newcommand{\di}{\partial}\newcommand{\bl}{\langle}\newcommand{\br}{\rangle}
\newcommand{\suml}{\sum\limits}\newcommand{\cAv}{A^\vee}
\newcommand{\quotient}[2]{{\left.\raisebox{0.4ex}{$#1$}\!\!\middle/\!\!\raisebox{-0.4ex}{$#2$}\right.}}

\newcommand{\tB}{\tilde{B}}\newcommand{\tC}{\tilde{C}}\newcommand{\tV}{\tilde{V}}

\newcommand{\ga}{\alpha}\newcommand{\gga}{\gamma}\newcommand{\gb}{\beta}\newcommand{\gd}{\delta}

\newcommand{\gl}{\lambda}\newcommand{\gL}{\Lambda}

\newcommand{\ber}{\begin{array}{l}}\newcommand{\eer}{\end{array}}
\newcommand{\bpm}{\begin{pmatrix}}\newcommand{\epm}{\end{pmatrix}}
\newcommand{\bM}{\begin{matrix}}\newcommand{\eM}{\end{matrix}}

\def\bull{\vrule height .9ex width .9ex depth -.1ex }

\newcommand{\beq}{\begin{equation}}\newcommand{\eeq}{\end{equation}}
\newtheorem{Lemma}{Lemma}[section]\newcommand{\bel}{\begin{Lemma}}\newcommand{\eel}{\end{Lemma}}
\newtheorem{Example}[Lemma]{Example}\newcommand{\bex}{\begin{Example}\rm}\newcommand{\eex}{\end{Example}}
\newtheorem{Proposition}[Lemma]{Proposition}\newcommand{\bprop}{\begin{Proposition}}\newcommand{\eprop}{\end{Proposition}}
\newcommand{\bpr}{~\\{\em Proof.~}}\newcommand{\epr}{$\bull$\\}
\newtheorem{Theorem}[Lemma]{Theorem}\newcommand{\bthe}{\begin{Theorem}}\newcommand{\ethe}{\end{Theorem}}
\newtheorem{Definition}[Lemma]{Definition}\newcommand{\bed}{\begin{Definition}}\newcommand{\eed}{\end{Definition}}
\newtheorem{Remark}[Lemma]{Remark}\newcommand{\beR}{\begin{Remark}\rm}\newcommand{\eeR}{\end{Remark}}
\newtheorem{Corollary}[Lemma]{Corollary}\newcommand{\bcor}{\begin{Corollary}}\newcommand{\ecor}{\end{Corollary}}
\newcommand{\bet}{\begin{tabular}{cccccccc}}\newcommand{\eet}{\end{tabular}}

\newcommand{\St}{{\rm St}}

\newcommand{\one}{{1\hspace{-0.1cm}\rm I}}\newcommand{\zero}{{0\put(-3.5,0){\line(0,1){7.5}}}}

\newcommand{\Ker}{{\rm Ker}}

\newcommand{\Ak}{{A^{(k)}}}
\newcommand{\Lie}{{\mathcal{L}ie}}
\newcommand{\wrt}{with respect to }

\font\frak=eufm10 scaled \magstep1
\def\frm{\mbox{\frak m}}

\title{N\MakeLowercase{ormal forms of matrices over the ring of formal series}}
\author[G.B\MakeLowercase{elitskii}]{G\MakeLowercase{enrich} B\MakeLowercase{elitskii}}
\address{Department of Mathematics, Ben Gurion University of the Negev, P.O.B. 653, Be'er Sheva 84105, Israel.}
\email{genrich@math.bgu.ac.il}
\author[D.K\MakeLowercase{erner}]{D\MakeLowercase{mitry} K\MakeLowercase{erner}}
\address{Department of Mathematics, Ben Gurion University of the Negev, P.O.B. 653, Be'er Sheva 84105, Israel.}
\email{kernerdm@math.bgu.ac.il}
\date{\today}
\thanks{D.K. was supported by the Skirball postdoctoral fellowship of the Center of Advanced Studies
in Mathematics (Mathematics Department of Ben Gurion University,
Israel). }

\begin{document}\maketitle\setcounter{secnumdepth}{6} \setcounter{tocdepth}{1}
\begin{abstract}
Matrices over the ring of formal power series are considered. Normal forms \wrt  various sub-groups of the two-sided transformations are constructed.
The construction is based on the special property of the action: it
induces a filtration  by projectors on sub-spaces of polynomial maps.
\end{abstract}
\tableofcontents

\section{Introduction}
Let $ G $ be a group acting on a space $X$. Then $X$ fibres into the disjoint union of $G$-orbits:
\beq
X = \bigcup\limits_{x\in X} Gx, \hspace{1cm} Gx = \{(gx) \in X |   g \in G\}
\eeq
Recall that a subset $ N \subset X $ is called a {\em normal form \wrt  the action of $ G $}
if $ N $ intersects all the orbits. An element $ z \in N \bigcap Gx $ is called a {\em normal form of $ x. $}
If the intersection is one point set for every $ x, $ then the normal form is called {\em canonical}.
A canonical normal form
solves the problem of classification of orbits: elements $ x_1 $ and $ x_2 $ lie in the same orbit if and only if
their canonical forms coincide.
The are many examples of normal and canonical forms in linear algebra \cite{Gantmacher-book}
and in analysis \cite{AGLV}.

Here we consider normal forms in spaces of matrices over the ring of formal power series.
Various groups of invertible matrices act on these spaces. Our considerations are based on the existence
of filtration preserved by these natural actions.
\subsection{Matrices of formal series}
Let $\K$ be a field of zero characteristic., not necessarily algebraically closed, \mbox{char$\K=0$}. Fix some natural $m,n,p$, denote by
\beq
Mat(m,n,p)= {\rm Mat}(m\times n,\ \K[[x_1 \ldots x_p]])
\eeq
the space of all $m\times n$ {\em formal matrices}, i.e. matrices $A(x)=(a_{ij}(x))$ whose
entries \mbox{$a_{ij}\in\K[[x_1 \ldots x_p]]$} are formal series of $ p $ variables over $ \K. $

Every matrix $ A \in Mat(m,n,p) $ can be represented as the formal series
$A(x)=\suml_{|I|=0}^\infty A_Ix^I$.
Here $I=(I_1 \ldots I_p) $ is an integer multi-index,
\beq
I_j \geq 0, \hspace{1cm} |I| = I_1 + \ldots + I_p, \hspace{1cm} x^I = x_1^{I_1} \ldots x_p^{I_p},
\eeq
and $ A_I $ is an $ m\times n $ matrix over $ \K. $

The space $ Mat(m,m,p)$ is a $\K$-algebra. Its element
$U(x) = U_0 + \suml_{|I|\geq 1} U_ix^I$ is invertible if and only if det$U_0\neq0$.
The subset $ GL(m,\K[[x_1 \ldots x_p]]) $ of all invertible elements is a group acting on the
space $ Mat(m,n,p) $ by multiplication from the left. Similarly, the
group $ GL(n,\K[[x_1 \ldots x_p]]) $ acts by multiplication from the right. The direct product
\beq
G(m,n,p)=GL(m,\K[[x_1 \ldots x_p]])\times GL(n,\K[[x_1 \ldots x_p]])
\eeq
acts from the two sides:
\beq
(g.A)(x) = U(x)A(x)V^{-1}(x), \hspace{1cm} g = (U,V).
\eeq
\subsection{Types of equivalence}
Formal matrices $A(x)$ and $B(x)$ are {\em two-sided equivalent} if they are equivalent by the
 action of the group $G(m,n,p)$. Similarly, they are {\em left equivalent} if
\beq
B(x) = U(x)A(x),\hspace{1cm} U \in GL(m,\K[[x_1 \ldots x_p]])
\eeq
and {\em right equivalent} if
\beq
B(x) = A(x)V^{-1}(x),\hspace{1cm} V \in GL(n,\K[[x_1 \ldots x_p]]).
\eeq
If the matrices are square, $ m = n, $ then they are {\em conjugate} if
\beq
B(x) = U(x)A(x)U^{-1}(x),\hspace{1cm} U \in GL(m,\K[[x_1 \ldots x_p]]),
\eeq
and are {\em congruent} if
\beq
B(x) = U(x)A(x)(U(x))^T,\hspace{1cm} U \in GL(m,\K[[x_1 \ldots x_p]]).
\eeq
Here $U^T$ means the transposition.
\\
\\
\\
The main aim of the paper is the construction of normal forms \wrt  the action of various
sub-groups $ G \subset G(m,n,p)$.
We extend here the approach suggested in \cite{Belitskii-1979-1}, \cite{Belitskii-1979-2}
 for locally analytic problems.
\subsection{One variable case}\label{Sec.One.Variable.Case}
In the simplest case of a single variable, i.e. $p=1$, a normal form \wrt  the two-sided
equivalence can be stated immediately. Set
\beq
N = \left\{\bpm R(x)&0\\0&0\epm\Big| \ \  R(x) = x^{k_1}\one_{m_1} \oplus \ldots\oplus x^{k_s}\one_{m_s}\right\}
\eeq
where $\one_r$ is the identity $ r\times r $ matrix, and
\beq
m_1 + \ldots + m_s \leq \min(m,n), \ k_s>k_{s-1}>\ldots>k_1.
\eeq
\bprop\cite{Birkhoff-1913}\cite{Grothendieck-1957}\label{Thm.Normal.Form.One.Variable.Case}
Every matrix of formal series in one variable is two-sided equivalent to a unique matrix in $N$.
\eprop
One can easily obtain similar normal form statements for other types of equivalences.
\subsection{Results}
Of course, the case $ p\geq 2 $ is essentially more complicated. Let us present some corollaries of our construction.

Let $\K\subset \C$. Given a polynomial $m\times n$ matrix
\beq
A(x) = \suml_{|I|\leq k} A_Ix^I, \hspace{1cm}    A_I\in{\rm Mat}(m\times n,\K),
\eeq
we introduce the differential operator
$A^*\left(\frac{\di}{\di x}\right):= \suml_{|I|\leq k}A_I^*\frac{\di^{|I|}}{\di x^I}$,
where $A^*_I:=\overline{A_I^T}$  is the conjugation.
It acts on the space $ Mat(m,n,p) $ by
\beq
A^*\left(\frac{\di}{\di x}\right)b(x)=\suml_{|I|\leq k}A_I^*\frac{\di^{|I|}b(x)}{\di x^I}.
\eeq

\bthe\label{Thm.Normal.Form.Main.Theorem}
 Let $\Ak(x)=\suml_{|I|=k}A_Ix^I, ~ A_I \in {\rm Mat}(m\times n,\K)$
be a homogeneous polynomial matrix of degree $k$, and let $A(x)=\Ak(x)+\big(\text{terms of orders}\geq k+1\big)$.
Then the matrix $A$ is two-sided equivalent to a matrix $B(x)=\Ak(x)+b(x)$ satisfying both of the relations
\beq\label{Eq.Relations.For.Normal.Form}
\Ak^*\left(\frac{\di}{\di x}\right)b(x) = 0,\hspace{1cm}
(\Ak^T)^*\left(\frac{\di}{\di x}\right) b^T(x) = 0
\eeq
\ethe
Hence, the set $N\subset Mat(m,n,p)$, consisting of matrices $B$ satisfying
equation (\ref{Eq.Relations.For.Normal.Form}), is a normal form \wrt  the action of the group $G(m,n,p)$.
It is not canonical: matrices $ B,\tB \in N $ may be two-sided equivalent.
\\
This theorem follows from general Theorem \ref{Thm.Normal.Form.from.Inner.Product} proved
in \S\ref{Sec.Normal.Form.from.Inner.Product}. It presents normal forms \wrt  various
sub-groups $ G \subset G(m,n,p)$, including
the left and the right equivalence, the conjugacy and so on.
 These normal forms turn out to be canonical \wrt  the ``unipotent parts'' of
 the respective groups consisting of transformations with identity linear part.
\\
\\
Two trivial cases:
\li
For $k=0$, the matrix $A^{(0)}(x)=A_0$ is a constant matrix, and we
obtain:
\bcor\label{Thm.Introduction.Corollary.Minimal.Resolution}
 Every matrix $A(x)=A_0+(\text{terms of orders}>0)$
 is two-sided equivalent to a matrix $ B(x)=A_0+b(x)$ satisfying $A_0^*b(x)=0$ and $b(x)A_0^*=0$.
\ecor
In commutative algebra the analogous statement is known as the reduction to the minimal resolution \cite{Eisenbud-book}.
Similar statements for other types of equivalence are obtained in \S\ref{Sec.Corollaries.And.Examples}.
\li For one variable case, i.e. $p=1$, Theorem \ref{Thm.Normal.Form.Main.Theorem} implies
Proposition \ref{Thm.Normal.Form.One.Variable.Case}
\\\\\\
Let now $m=n$ and $\Ak(x)={\rm diag}\big(l_1(x) \ldots l_m(x)\big)$ a diagonal matrix with homogeneous
polynomials $l_i(x)$ of degree $ k$.  If $b(x) = (b_{ij}(x))$, then
\beq
\Ak^*\left(\frac{\di}{\di x}\right)b(x) =\left\{ l_i^*\left(\frac{\di}{\di x}\right)b_{ij}(x)\right\}_{i,j=1}^m
 ~ \text{ and }~
(\Ak^T)^*\left(\frac{\di}{\di x}\right)b^T(x) =\left\{ l_j^*\left(\frac{\di}{\di x}\right)b_{ij}(x)\right\}_{i,j=1}^m
\eeq
Here
\beq
l^*_i\left(\frac{\di}{\di x}\right)=\suml_{|I|=k}\overline{(l_i)}_I\frac{\di^{|I|}}{\di x^I}.
\eeq
Then Theorem \ref{Thm.Normal.Form.Main.Theorem} leads to
\bcor\label{Thm.Introduction.Corollary.For.Diag.Matrix} Every matrix $A(x)=\Ak(x)+(\text{terms of orders}>k)$
 with $\Ak(x)={\rm diag}(l_1(x) \ldots \l_m(x)) $ is two-sided equivalent to a matrix
$ B(x)=\Ak(x)+(b_{ij}(x)) $ satisfying
\beq
l_i^*\left(\frac{\di}{\di x}\right)b_{ij}(x) = 0, \hspace{1cm}
l_j^*\left(\frac{\di}{\di x}\right)b_{ij}(x) = 0
\eeq
\ecor
\subsection{Matrices with formal/locally convergent/rational entries}
In many applications the matrix functions are considered as local objects, defined near the origin. (For example in
singularity theory or local algebraic geometry.) Accordingly one has various notions of locality:
\li the neighborhoods of Zariski topology (corresponding to the rational functions, regular at the origin),
\li the neighborhoods of the classical topology (locally convergent functions)
\li formal neighborhoods (formal series).

While the formal neighborhoods are better for theoretical considerations (e.g. no issues of convergence),
in practice one works usually with the locally convergent series or rational functions.
Correspondingly one has various comparison questions of the three cases. We discuss this
in \S\ref{Sec.Dependence.On.The.Base.Ring}.
\section{Filtration in the spaces of formal matrices.}\label{Sec.Filtration.in.the.space.of.formal.matrices}
Given $ A \in Mat(m,n,p) $ and $j\in\N$,  consider the {\em $j$'th jet} of the matrix,
$\pi_jA:= \suml_{|I|\leq j}A_Ix^I$. Then $\{\pi_j\}$ is an increasing system of projectors, i.e.
$\pi_i\pi_j = \pi_j\pi_i = \pi_i, \ i \leq j$.
The image $Mat_j(m,n,p)$ of the projector $\pi_j$ consists of all polynomial matrices of degree $\leq j$.
\\
\\
\\
The {\em homogeneous summand} $A^{(j)}(x):= \suml_{|I|=j}A_Ix^I$ can be represented in the form
\beq
A^{(j)} = (\pi_j - \pi_{j-1})A, \pi_{-1} = 0.
\eeq
Correspondingly, the projector to the j'th homogeneous component is $\pi^{(j)}:=\pi_j-\pi_{j-1}$.
 Its image, $Mat^{(j)}(m,n,p)$,  consists of all the homogeneous matrices of degree $j$. Thus
\beq
\pi^{(j)}\pi^{(i)}=\Big\{\bM \zero,& i\neq j\\\pi^{(i)},&j=i \eM
\eeq
\\
\\
\\
The system $ \{\pi_j\} $ generates a ``sequential topology'': a sequence
$ \{ A_k\} $ {\em converges} to $ A $ if
\beq
\forall j\in\mN \text{ there exists } k_0(j)\in\mN\text{ such that }\forall k \geq k_0(j):\  \pi_jA_k = \pi_jA,
\eeq

The ``convergence Cauchy criterion'' states that a sequence $\{A_k\}$ converges if and only if it stabilizes:
\beq
\pi_jA_k = \pi_jA_{k'}, k,k' \geq k_0(j), j = 0,1, \ldots .
\eeq
Then the matrix $ A $ with the jets
\beq
\pi_jA = \pi_jA_k, \ k \geq k_0(j)
\eeq
is the limit.

A subset $S\subset Mat(m,n,p)$ is {\em closed} if the limit of every converging sequence
$\{A_k\}\subset S $ belongs to $S$.

\bex
$\bullet$ A one-point set is a closed subset in $ Mat(m,n,p), $ and is not open:
for every $ A \in Mat(m,n,p) $ the complement
$ Mat(m,n,p)\setminus\{ A\} $ is not closed.

\li Let $ P \in Mat(m,n,p) $ and let $ s \geq 0 $ be a fixed integer. The subset
\beq
\{ A \in Mat(m,n,p) ~ | ~ \pi_sA =\pi_sP\}
\eeq
is simultaneously closed and open.
\eex~
\\
\\
\\
Consider now the direct product $Mat(m,m,p) \times Mat(n,n,p)$.
Each $A\in Mat(m,m,p) \times Mat(n,n,p)$ determines the two-sided action, so we will write
\beq
A = (A_l,A_r), ~ A_l \in Mat(m,m,p), ~ A_r \in Mat(n,n,p).
\eeq
The j-jet projectors to the subspaces of polynomials, $\{\pi_j\}$, act on  $Mat(m,m,p) \times Mat(n,n,p)$ as previously.
In what follows we denote the projectors by the same letters: \mbox{$\pi_jA = (\pi_jA_l,\pi_jA_r)$.}
\section{Lie groups of formal transformations.}\label{Sec.Lie.Groups}
The group $ GL(m,\K[[x_1 \ldots x_p]]) $ is a ``countably dimensional Lie group'', possessing an exponential map,
as we explain now.
\subsection{The exponential and logarithmic maps.}
Let $\gl\in Mat(m,n,p)$ be constant term free, i.e. $ \pi_0\gl = 0$. Then the sequence
\beq
U_k(x)=\suml_{i=0}^k\frac{\gl^i(x)}{i!}
\eeq
converges to a matrix $ \exp\gl \in GL(m,\K[[x_1 \ldots x_p]]) $ such that $ \pi_0\exp\gl = $id.

The image of the map $\exp$ coincides with the sub-group of all elements $ U \in GL(m,\K[[x_1 \ldots x_p]]) $
such that $\pi_0U=\one$. Indeed, let $U(x)=\one + u(x)$ with $\pi_0u=\zero$. Then the sequence
\beq
\gl_k(x) = \suml_{i=1}^k \frac{(-1)^{i-1}}{i}u^i(x)
\eeq
converges to the matrix $ \log U \in Mat(m,m,p)$ with $\pi_0\log U =\zero$, satisfying $\exp(\log U) = U$.
\subsection{Lie groups and their algebras.}\label{Sec.Lie.Groups.and.their.Algebras}
Given an $m\times m$ matrix $\gl$ with $\pi_0\gl=\zero$, the map
\beq
\Phi : \K \rightarrow GL(m,\K[[x_1 \ldots x_p]])
\eeq
defined by $\Phi(t) = \exp(t\gl)$ is a {\em flow}, i.e. satisfies the relation
$\Phi(t_1 + t_2) = \Phi(t_1)\Phi(t_2)$.
\\
Its image is a one-parametric sub-group. In this sense $ GL(m,\K[[x_1 \ldots x_p]]) $ is
a Lie group. The Lie algebra of this group is
\beq
\Lie:= \{ \gl \in Mat(m,m,p) : \pi_0\gl =\zero \}
\eeq
with the usual commutator $[\gl_1,\gl_2]:=\gl_1\gl_2 - \gl_2\gl_1.$

Similarly, the group $G(m,n,p)=GL(m,\K[[x_1 \ldots x_p]])\times GL(n,\K[[x_1 \ldots x_p]])$
is Lie with the Lie algebra
\beq
\Lie= \{ \gl = (\gl_l,\gl_r) \in Mat(m,m,p)\times Mat(n,n,p), \ \ \pi_0\gl_l =\zero,\  \pi_0\gl_r =\zero \}
\eeq
and the commutator $[\nu,\gl]:=([\nu_l,\gl_l], [\nu_r,\gl_r]).$
\\
\\
\\
More generally, let $ G \subset G(m,n,p) $ be a sub-group. Define its subgroup of matrices
that are "close to being idempotent":
\beq
G^0 = \{ g = (U,V) \in G| \ \pi_0U =\one,\  \pi_0V =\one\}.
\eeq
Then $ G^0 $ is contained in the image of the exponential map.
\bed
A closed sub-group $ G \subset G(m,n,p) $ is called Lie if for
every $ g = \exp\gl \in G^0 $ the one-parametric group $\{\exp(t\gl),~ t\in\K\}$ is contained in $ G $
entirely.
\eed
In the sequential topology the sub-group $ G^0 $ is an open neighborhood of the identity in the Lie group $G$.

Note, if $G$ is a Lie group then so is $G^0$ and their Lie algebras coincide.

Given a Lie group $ G $ we denote by $ \Lie(G) $ its Lie algebra. If $ \gl=(\gl_l,\gl_r) \in
\Lie(G), $ then $g=\exp\gl \in G^0 $ and $g.A = \exp\gl_lA\exp(-\gl_r)$.

\bex
Here are some commonly used Lie groups and their algebras.
\li We denote the ``biggest'' group $G(m,n,p)$  by $G_{lr}$.
Matrices equivalent \wrt  this group are called {\em two-sided equivalent}. Obviously,
\beq
\Lie(G_{lr}) =\{ \gl \in Mat(m,m,p)\times Mat(n,n,p):\pi_0\gl_l=\zero=\pi\gl_r \}.
\eeq
$\bullet$ The group $G_l$ of the left transformations $g=(U,\one)$ is
Lie too. We call $G_l$-equivalent matrices {\em left equivalent}. The corresponding Lie algebra is
\beq
\Lie(G_l)=\{\gl : \gl_r =\zero,\  \pi_0\gl_l =\zero \}
\eeq
$\bullet$ The group $G_r$ of the right transformations $ g = (\one,V) $ is Lie, and its Lie algebra is
\beq
\Lie(G_r) = \{ \gl : \gl_l=\zero,\  \pi_0\gl_r=\zero \}.
\eeq
$\bullet$ Matrices equivalent \wrt  the diagonal subgroup
\beq
G_c = \{ g = (U,U^{-1}) \}\subset G(m,m,p)\times G(m,m,p)
\eeq
are called {\em conjugate}. Obviously, $\Lie(G_c) = \{ \nu : \gl_l=-\gl_r,\  \pi_0\gl_l=\zero\}$.
\li  The group $G_T = \{g=(U,U^T)\}\subset G(m,m,p)\times G(m,m,p)$  is Lie too, and
\beq
\Lie(G_T) =\{ \gl : \gl_r=\gl_l^T, \pi_0\gl_l=\zero\}.
\eeq
$G_T$-equivalent matrices are called
{\em congruent}.
\eex
~\\
The action of any such group $G\subset G(m,n,p) $ on $Mat(m,n,p) $ is consistent with the projectors:
\beq\label{Eq.Group.Action.Consistent.WithProjectors}
\pi_j(g.A) = \pi_jg.\pi_jA, ~ j = 0,1, \ldots
\eeq
If $ g \in G^0, $ then, in addition, the implication
\beq\label{Eq.From(j-1)jet.to.j.jet}
\pi_{j-1}A = \pi_{j-1}B\ \Rightarrow\\\pi_j(g.A - g.B) = \pi_j(A - B)
\eeq
holds. These properties allow us to adjust jets of a given formal matrix $ A$ successively.
\subsection{Linearization of the Lie groups actions.}
First we study the {\em stabilizer} group of a given matrix, $\St(A):=\{g| ~ g.A = A \}\subset G^0(m,n,p)$.
\bprop
The stabilizer is a Lie group.
\eprop
\bpr Let $g=\exp\nu\in\St(A)$. Consider the function
\beq\ber
f : \K \rightarrow Mat(m,n,p)
\\f(t) = \exp(t\nu).A - A
\eer\eeq
Since $ \St(A) $ is a sub-group, $ f(n) = 0 $ for $ n = 0,\pm 1, \ldots . $ Hence, the functions
\beq
f_j(t) = \pi_jf(t)
\eeq
vanish on an infinite set of points. Being polynomials, they are zero identically.
Hence, $ \exp t\nu \in \St(A) $ for all $ t \in \K. $
\epr

Similarly, the sub-groups of matrices stabilizing specific jets:
\beq
\St_j(A):= \{ g \in G^0 ~ |~  \pi_jg.A = \pi_jA \}, j = 0,1, \ldots
\eeq
are Lie too.
\bed\label{Def.Linear.Map.On.Lie.Algebra}
Given $ \nu \in \Lie(G_{lr}), $ define the linear map
$S_A'(0)\nu:= f'(0) = \frac{d}{dt}(\exp t\nu.A - A)|_{t=0}$:
\beq
\bM S_A'(0) : {\rm Lie}(G_{lr}) \rightarrow Mat(m,n,p)\\\nu=(\nu_l,\nu_r)\mapsto \nu A=\nu_lA - A\nu_r\eM
\eeq
It is called the linearization of the action of $G_{lr}$.
\eed
\bel\label{Thm.Stj(A).vs.Ker(piS'A)}
1. An element $\exp\nu$ belongs to $ \St_j(A) $ if and only if $ \nu \in \Ker\Big( \pi_jS_A'(0)\Big)$,
i.e. $\pi_j(\nu A)=\zero$.
\\2. An element $\exp\nu$ belongs to $\St(A) $ if and only if $\nu\in\Ker S_A'(0)$, i.e. $\nu A=\zero$.
\\3. If $ \exp\nu \in \St_{j-1}(A), $ then
\beq
\pi_j\exp\nu.A = \pi_jA + \pi_j(S_A'(0)\nu)
\eeq

\eel
\bpr 1. If $\exp\nu\in St_j(A)$ then $\pi_j(\exp t\nu.A-A)=\zero$ for all $t\in \K$. Thus $ \pi_jS_A'(0)\nu=\zero$.

Conversely, let the latter equality hold. Then $\pi_j(\nu^kA)=\zero$, i.e. $\pi_j\nu_l^kA=\pi_jA\nu_r^k$, for any $k$.
 Hence, $\pi_j\exp\nu_l.A = \pi_jA\exp\nu_r$, meaning that $ \exp\nu \in \St_j(A)$.
\\
\\
2. If $\exp\nu\in St(A)$ then $\exp\nu \in \St_j(A)$, hence $\pi_jS_A'(0)\nu = 0$ for all $ j = 0,1, \ldots$.
\\
Conversely, let $ S_A'(0)\nu = 0. $ Then $ \exp\nu \in \St_j(A) $ for all $ j, $ implying $ \exp\nu \in St(A). $
\\\\
3. Let $ \exp\nu \in \St_{j-1}(A). $ Then $ \exp t\nu \in \St_{j-1}(A) $ for all $ t \in \K. $
It follows from equation (\ref{Eq.From(j-1)jet.to.j.jet}) that
\beq
\pi_j(\exp t\nu.\exp s\nu.A-\exp t\nu.A) =\pi_j(\exp s\nu.A - A)
\eeq
for all $ t,s \in \K. $ Hence,
\beq
f_j(t+s)=\pi_j(\exp(t+s)\nu.A-A)=\pi_j(\exp t\nu.\exp\nu s.A-\exp t\nu.A) + \pi_j(\exp t\nu.A - A)=f_j(t)+f_j(s),
\eeq
meaning that $ f_j $ is an additive function. Being a polynomial (and as $char(\Bbb K)=0$) it is linear:
\beq
f_j(t) = f_j'(0)t = \pi_jS_A'(0)t\nu
\eeq
\epr

\section{Determinacy by jets.}\label{Sec.Determinacy.By.Jets}
\subsection{Jet-by-jet equivalence.}
Let $ G \subset G(m,n,p) $ be a sub-group. Formal matrices $A(x) $ and $B(x) $ are called
{\em jet-by-jet $ G$-equivalent} if there is a sequence $ \{ g_j\} \subset G $ such that
\beq\label{Eq.Jet.by.jet.Equivalence}
\pi_jg_j.A = \pi_jB, \ j = 0,1, \ldots .
\eeq

A subgroup $G\subset G(m,n,p)$ is called {\em countably algebraic} over $\K$, if it is determined by a finite
or infinite system of polynomial equations in the matrix coefficients of monomials.
More precisely, this means:

$g=(U,V)\in G$, with $U(x) = \suml_I U_Ix^I$ and $V(x)=\suml_I V_Ix^I$, if and only if
\beq
P_k(\{ U_I\},\{ V_I\}) = 0, |I| \leq N_k,\  k = 1,2, \ldots
\eeq
for some polynomials $ P_k $ over $\K$ and numbers $N_k$.

Obviously, the groups $G_l,G_r,G_{lr},G_c,G_T$ are countably algebraic.
\bthe\label{Thm.Jet.by.Jet.equivalence.implies.Equivalence} Let the field $\K$ be algebraically closed, and let $G$
be a countably algebraic group.
Then jet-by-jet $G$-equivalence implies $G$-equivalence.
\ethe
\bpr Since $G$ is countably algebraic, the condition (\ref{Eq.Jet.by.jet.Equivalence}) is an infinite system of polynomial
equations \wrt  coefficients of the transformations $g_j$. By Lang's Theorem there is a
 common solution $g$ satisfying all equations.
\epr

Let us remind the {\em Lang Theorem} \cite{Lang-1952}.
\bthe Consider an infinite polynomial system
\beq\label{Eq.Langs.Theorem}
h_j(a_1 \ldots a_{m_j})= 0, ~ j = 0,1, \ldots, \hspace{1cm} m_j \rightarrow \infty
\eeq
over a field $ \K. $ Assume that for every $ k = 0,1,2, \ldots $ the finite sub-system
\beq
h_j(a_1 \ldots a_{m_j}) = 0, \ j = 0,1, \ldots, k
\eeq
is solvable. If $\K$ is algebraically closed, then the total initial system
has a solution $ \overline{a}=(a_1,a_2, \ldots ) $.
\ethe
If $\K$ is not closed the statement fails.
\bex Let $\K=\R$. Consider the system of equations:
\beq\label{Eq.Infinite.Syst.Eqs}
a_j^2 = a_1 - j, \ j = 1,2, \ldots .
\eeq
Given an integer $k\geq 0$, the sequence $a_1 = k$, $a_j=\pm\sqrt{k-j}$ is a solution of the system for $j\leq k$.
However, the infinite system (\ref{Eq.Infinite.Syst.Eqs}) has no real solution $ a_1,a_2, \ldots$. On the other hand,
for every $ a_1 \in \C $ the sequence
\beq
a_j = \pm\sqrt{a_1-j}, \ \ \ j = 2,3, \ldots
\eeq
is a complex-valued solution of (\ref{Eq.Infinite.Syst.Eqs}).
\eex
~\\\\
However, at least for some groups the statement of the last theorem is valid over an arbitrary field.
\bprop
Let $G$ be one of the groups $G_{lr}$, $G_l$, $G_r$, $G_c$.
Then the jet-by-jet $G$-equivalence implies $G$-equivalence.
\eprop
\bpr
Let the matrices $A(x)$ and $B(x)$ be jet-by-jet $G$-equivalent over $\K$.
 Let $\K\subset\bK$ be the algebraic closure.
Then, by Theorem \ref{Thm.Jet.by.Jet.equivalence.implies.Equivalence}, the matrices are $G$-equivalent
over $\bK$, i.e. $U(x)A(x) = B(x)V(x)$ with $U(0)$, $V(0)$ non-degenerate.

Let $\{w_\ga\}_{\ga}$ be a Hamel basis of $\bK$ as a vector space over $\K$, i.e. a maximal set
of $\K$-linearly independent elements, cf. \cite[pg.53]{Rudin-book}.
So, any element of $\bK$ is presentable as $\sum a_\ga w_\ga$, for $a_\ga\in\K$
and the sum is {\em finite}. Thus any series $f\in\bK[[x]]$ decomposes:
\beq
f=\suml_{j=0}^\infty\suml_{deg(I)=j}\suml_\ga a_{I,\ga} w_\ga x^I
\eeq
 Here $\suml_{j=0}^\infty\suml_{deg(I)=j}a_{I,\ga} x^I\in\K[[x]]$. Note that for each fixed $I$
the inner sum $\suml_\ga a_{I,\ga} w_\ga x^I$ is finite, hence there is no problem of convergence.

Similarly decomposes every matrix with entries in $\bK[[x]]$.
 Let $\{U_\ga(x)\}$ and $\{V_\ga(x)\}$ be the projections of $U(x)$,$V(x)$ onto the ''$w_\ga\K[[x]]$'' subspaces,
 i.e. matrices with entries in $\K[[x]]$.
Hence from $U(x)A(x) = B(x)V(x)$ one has: $U_\ga(x)A(x)=B(x)V_\ga(x)$ for any $\ga$.

We claim that there exists a sequence of numbers $\{\gl_\ga\in\K\}$ such that
$\sum \gl_\ga U_\ga(0)$ and $\sum \gl_\ga V_\ga(0)$ are non-degenerate matrices.
Then $\big(\sum \gl_\ga U_\ga(x)\big) A(x)=B(x)\big(\sum \gl_\ga V_\ga(0)\big)$ proves the statement.

Indeed, consider the polynomials $\det(\sum y_\ga U_\ga(0))$ and $\det(\sum y_\ga V_\ga(0))$ where $\{y_\ga\}$ are
independent variables. (As previously there is a finite number of variables.)
As the matrices $U(0)=\sum w_\ga U_\ga(0)$ and $V(0)=\sum w_\ga V_\ga(0)$ are non-degenerate these polynomial are
not identically zero.
Thus they are not identically zero for some value $y_1=\gl_1\in\K$. Fix this value then there exists $y_2=\gl_2\in\K$
such that $\det\big(\sum y_\ga U_\ga(0)\big)|_{\substack{y_1=\gl_1\\y_2=\gl_2}}\not\equiv0$ and
$\det\big(\sum y_\ga V_\ga(0)\big)|_{\substack{y_1=\gl_1\\y_2=\gl_2}}\not\equiv0$.
Continue by induction to build the needed (finite) sequence $\{\gl_\ga\in\K\}$.
\epr

Obviously, $G_T$-equivalence of real matrices over $ \C $ does not imply their congruence over $ \R. $
Nevertheless, the statement of Theorem \ref{Thm.Jet.by.Jet.equivalence.implies.Equivalence} is true over an
arbitrary field at least for the unipotent part $ G^0 $ of a Lie group $ G. $
\bthe\label{Thm.Jet.by.Jet.equivalence.implies.Equivalence.For.$G^0$} Let $G$ be a Lie group over an arbitrary $\K$.
Then the jet-by-jet $ G^0$-equivalence implies $ G^0$-equivalence.
\ethe
\bpr
Let a matrix $ B $ be jet-by-jet $ G^0$-equivalent to $A$. Then the sets
\beq
M_k=\{\nu\in\Lie(G) ~ | ~ \pi_k\exp\nu.A =\pi_kB\}
\eeq
are non-empty and decrease: $ M_{k+1} \subset M_k$.
\\{\bf Step 1.} We prove that the sequence of jets stabilizes:
\beq
\pi_jM_{k+1} = \pi_jM_k, \ k \geq k_0(j),\  j = 0,1, \ldots .
\eeq

This is immediate if $\K$ is algebraically closed.
Indeed, for any finite $j$ the subset $\pi_j M_k\subset \pi_j Mat(n,m,p)$
is an algebraic subvariety of finite dimension.
Hence the decreasing sequence $\pi_jM_{k+1}\subset\pi_jM_k$ necessarily stabilizes.
Namely, for any $j$: $\pi_jM_{k+1} = \pi_jM_k$, for $k\geq k_0(j)$.
\\
\\~

For an arbitrary field we prove as follows.
Fix a sequence $\{\nu_k \in M_k\}_{k =0,1, \ldots}$. Then for any $\nu\in M_k$ have:
$\exp(-\nu_k)\exp\nu\in St_k(A)$, i.e. there exists $\tau \in \Ker \pi_kS_A'(0)$ such that
$\exp\nu = \exp\nu_k\exp\tau$. And conversely, if $\exp\tau\in\St_k(A)$ then
$\exp\nu_k\exp\tau\in M_k$. This defines the set-theoretic bijection:
\beq\ber
\nu: \ \Ker \pi_kS_A'(0)\rightarrow M_k\\\tau\mapsto\exp(\nu_k)\exp(\tau)
\eer\eeq
Similarly the map $\nu_j:\pi_j\big(\Ker \pi_kS_A'(0)\big)\rightarrow \pi_j(M_k)$ is bijection too.

Now we have a decreasing sequence of vector spaces:
\beq
\pi_j(\L ie)\supset..\supset \pi_j\Big(Ker\pi_kS_A'(0)\Big)\supset\pi_j\Big(Ker\pi_{k+1}S_A'(0)\Big)\supset...
\eeq
They are subspaces of a finite dimensional space, hence the sequence stabilizes:
\beq
\cap_k\pi_j\Big(Ker\pi_kS_A'(0)\Big)=\pi_j\Big(Ker\pi_{k_0(j)}S_A'(0)\Big)
\eeq
By the bijection above we get $\cap_k\pi_jM_k=\pi_jM_{k_0(j)}$.
\\\\~

{\bf Step 2.}
Denote
\beq\label{Eq.S_j}
S_j = \pi_jM_k, \ k \geq N(j).
\eeq
We can assume the sequence $ N(j) $ is increasing. Then for $i\leq j$ :
\beq
\pi_iS_j = \pi_i\pi_jM_k =\pi_iM_k = S_i, \ k \geq N(j).
\eeq
i.e. the projection $S_j\to S_i$ is surjective.
Hence, we can choose successively elements $\nu_0=0$, $\nu_1\in S_1$, $\nu_2\in S_2$,  with $\pi_1\nu_2=\nu_1$,
then $\nu_3\in S_3$ with $\pi_2\nu_3=\nu_2$ and so on. The sequence $\{\nu_j\}$ converges. Its limit $\nu$ satisfies
\beq
\pi_k\exp\nu.A = \pi_kB, \ k = 0,1, \ldots .
\eeq
Therefore, $ A $ and $ B $ are $ G^0$-equivalent.
\epr

\subsection{Finite determinacy.}
Now we discuss when a finite jet determines the $G$-equivalence class of the matrix.
\\
\\
Let us recall, that a matrix $A(x) $ is called {\em $ k$-determined} \wrt  a group
 $ G $ if every matrix $B$ whose $k$-jet $\pi_k(B)$ equals $\pi_kA$ is $ G$-equivalent to $A$. The minimal such $k$
 is called the {\em order of determinacy} \wrt  the group $G$.
 A matrix is called {\em finitely determined} if it is $k$-determined with $k<\infty$.
 Otherwise the matrix is called {\em infinitely determined}.
\bprop
1. The order of determinacy (finite or infinite, \wrt any group) is invariant \wrt $G_{lr}$ action.
\\2. If a matrix is $k$-determined \wrt  a group $G$, then it is $G$-equivalent to a matrix whose entries
are polynomials of degrees at most $k$.
\eprop
\bpr
1. Suppose $A(x)$ is $k$-determined \wrt $G$. Let $B(x)_{\ge k}$ be a matrix with $\pi_{k-1}(B(x)_{\ge k})=\zero$.
 So the matrix $A(x)+B(x)_{\ge k}$ is $G$-equivalent to $A(x)$. Then for any $g\in G_{lr}$
the matrix $gA(x)$ is $k$-determined too:
\beq
gA(x)+B(x)_{\ge k}=g\big(A(x)+g^{-1}B(x)_{\ge k}\big)\stackrel{G}{\sim}gA(x)
\eeq
because $\pi_{k-1}\big(g^{-1}B(x)_{\ge k}\big)=\pi_{k-1}g^{-1}\pi_{k-1}\big(B(x)_{\ge k}\big)=\zero$.
\\2. Immediately.
\epr
\bex
Many matrices  are not $G_{lr}$ equivalent to polynomial matrices. Consider a $1\times 1$ matrix
$A=\{y-xf(x)\}$, where $f(x)$ is a locally analytic but not rational function.
For example $f(x)=\exp(x)$. Then the curve $\{y-xf(x)=0\}\subset\C^2$ is locally analytic but not algebraic.
Hence any equivalent matrix cannot be a polynomial (as it must define the same non-algebraic curve).
\eex

\bthe\label{Thm.Finite.Determinacy.Equivalent.Condition} Let $G$ be a Lie group over a field $\K$.
For $A\in Mat(m,n,p)$ let $S_A'(0)$ be the map as defined in \ref{Def.Linear.Map.On.Lie.Algebra}.
\\1. Assume that
$Mat^{(j)}(m,n,p)\subset \pi_jS_A'(0)(\Lie(G))=\pi_j\big(\Lie(G)A\big)$ for $j\ge k+1$.
Then the matrix $ A $ is $ k$-determined \wrt  $G^0$ and, as a consequence also
\wrt  $G$.
\\2. Conversely, if $A$ is $k$-determined \wrt  $G^0$ then $Mat^{(j)}(m,n,p)\subset \pi_jS_A'(0)(\Lie(G))$
for $j\ge k+1$.
\ethe
\bpr 1. Suppose the condition holds, let $B(x) = A(x) + P(x)$ and $\pi_kP(x)=0$.
We should show that $B$ is $G^0$-equivalent to $A$.

Let $ P^{(k+1)} = \pi_{k+1}S_A'(0)\nu_1$, for $\nu_1 \in \Lie(G)$. Set $ g_1 = \exp\nu_1. $ Then
\beq
\pi_{k+1}(g_1.A - B) =\pi_{k+1}(g_1A - A - P^{(k+1)}) = 0.
\eeq
Hence, $ g_1^{-1}B = A + P_1$ with $\pi_{k+1}P_1=0$.
Further, let $P_1^{(k+2)}=\pi_{k+2}S_A'(0)\nu_2$ for $\nu_2 \in \Lie(G)$.
Then $g_2^{-1}g_1^{-1}B = A + P_2$ with $\pi_{k+2}P_2=0$.

In general, for every $j=k+1,k+2\ldots$, there is a transformation $ h_j = g_j^{-1}g_{j-1}^{-1} \ldots
g_1^{-1} \in G^0 $ such that $ \pi_jh_jB = \pi_jA. $ This means that $ A $ and $ B $ are jet-by-jet $ G^0$-equivalent.
By Theorem 4.2 they are $ G^0$-equivalent.
\\
\\
2. Conversely, let $A$ be $ k$-determined \wrt  $G^0$. Then every matrix
$B(x)=A(x)+P^{(j)}(x)$, for $j\geq k+1$ is $G^0$-equivalent to $A$. Hence,
\beq
\exp\nu.A=A+P^{(j)}(x), \ \nu\in\Lie(G).
\eeq
Since $ \pi_{j-1}(A + P^{(j)}) =
\pi_{j-1}A, $ the element $ g = \exp\nu $ lies in the stabilizer $ \St_{j-1}(A)$.
By Lemma \ref{Thm.Stj(A).vs.Ker(piS'A)}
\beq
P^{(j)} = \pi_j(\exp\nu.A - A) =\pi_jS_A'(0)\nu.
\eeq
\epr
\bcor\label{Thm.Corollary.Finite.Determinacy}
 For a given $A$, suppose the equation $P=S_A'(0)\nu$ has a solution $\nu\in\Lie(G)$ for every formal matrix $P$ with $\pi_kP=0$.
 Then $A$ is $k$-determined \wrt  $G^0$.
\ecor
Indeed, if $P\in S_A'(0)(\Lie(G))$ for every matrix $P$ with $\pi_kP=0$, then the condition of the theorem is satisfied.
\bex Let $A=$const be a constant $ m\times n $ matrix. It is $ 0$-determined \wrt  the group $ G_l^0 $ if and
only if it is invertible from the left, i.e. $ \Ker A = \{ 0\}, $ or rank$ A = n. $ Otherwise it is not finitely determined.
Similarly, $ A $ is $ 0$-determined \wrt  the group $ G_r^0 $ if and only if it is invertible from the right, i.e.
rank$ A = m. $ The matrix $ A $ is $ 0$-determined \wrt  two-sided transformations if and only if rank$ A = \min(m,n). $
This means that the linear transformation given by $A$ is either surjective or injective.
Indeed, assume that the equation
\beq
\gl_lA + A\gl_r = P
\eeq
has a solution with an $m\times m$ matrix $\gl_l$ and an $n\times n$ matrix $\gl_r. $ Assume $ Ax = 0$, $x\neq 0$.
Then the equation takes on the form $ A\gl_rx = Px. $ Since the vector $ Px $ is arbitrary, rank$ A = m. $

However, a constant matrix is not finitely determined \wrt  conjugacy. Indeed, the equation
\beq
\gl A - A\gl = P
\eeq
is solvable only for $P$ with trace$(P)=0$.
\eex

\section{Normal forms}\label{Sec.Normal.Forms}
\subsection{Construction of the normal form}\label{Sec.Normal.Forms.Existence}
In this section we give a constructive description of the normal form. Or, in elementary terms: given $A\in Mat(m,n,p)$
and a group $G\subset G(m,n,p)$, how to reduce $A$ modulo the orbit $GA$.

As $Mat(m,n,p)$ is a vector space graded by the total degree, i.e. $Mat(m,n,p)=\oplus_j Mat(m,n,p)^{(j)}$,
it is natural to apply the jet-by-jet reduction. Namely, at the j'th step we adjust the j'th jet, preserving
the (j-1)'st jet.

Define the stabilizer
\beq
St_{j-1}(A):=\{g\in G| \ \pi_{j-1}(gA)=\pi_{j-1}(A)\}=\{\nu\in \Lie(G)| \ \pi_{j-1}(\nu A)=0\}
\eeq
Here the last equality is due to Lemma \ref{Thm.Stj(A).vs.Ker(piS'A)}.

The subgroup $St_{j-1}(A)\subset G$ defines the orbit $St_{j-1}(A)A$.
Consider its j'th jet:
\beq
\pi_j\Big(St_{j-1}(A)A\Big)=\{\pi_j(gA)| \ \pi_{j-1}(gA)=\pi_{j-1}(A)\}=
\{\pi_j(A)+\pi_j(\nu A)| \ \pi_{j-1}(\nu A)=0\}
\eeq
Again the second equality is due to Lemma \ref{Thm.Stj(A).vs.Ker(piS'A)}.

This defines the vector space
\beq
Mat(m,n,p)^{(j)}\supset V^{(j)}(A):=\{\pi_j(\nu A)| \ \nu\in \Lie(G), \ \pi_{j-1}(\nu A)=0\}
\eeq
Note that $V^{(j)}(A)=V^{(j)}(\pi_{j-1}A)$.
Let $W^{(j)}(A)$ be a complementary subspace, i.e. $ V^{(j)}(A)\oplus W^{(j)}(A)=Mat(m,n,p)^{(j)}$.
Define
\beq
N(G) = \{ B \in Mat(m,n,p) : B^{(j)}\in W^{(j)}(B), j=0,1,\ldots \}
\eeq
\bthe\label{Thm.Normal.Form.Exists}
1. The set $N(G)$ is a normal form \wrt  the action of $G$ on $Mat(m,n,p)$.
\\2. It is a canonical form \wrt  $ G^0$.
\ethe
\bpr
1. Given a matrix $A(x)$, we construct inductively $g\in G^0$ and $B\in Mat(m,n,p)$, such that $B=gA\in N(G)$.
Set $\pi_0B = A^{(0)}$. Suppose we have built $ g_{i-1} \in G^0 $ such that the jet
$B_{j-1}:=\pi_{j-1}(B)=\pi_{j-1}(g_{j-1}A)$ is in the normal form, i.e.
$ B_{j-1}^{(i)} \in W^{(i)}(B_{i-1})$, for $i \leq j-1$.
\\
\\~

Let $\pi_j(g_{j-1}A)=B_{j-1}+v_j+w_j$, where $v_j\in V^{(j)}(g_{j-1}A)$ and $w_j\in W^{(j)}(g_{j-1}A)$.
By construction $v_j=\pi_j(\nu_j A)$ for some $\nu_j$. Hence
\beq
\pi_j(\exp(-\nu_j)g_{j-1}A)=\pi_j((1-\nu_j)g_{j-1}A)=B_{j-1}+w_j=:B_j
\eeq
Thus define $g_j:=\exp(-\nu_j)g_{j-1}$ and we get that $B_j=\pi_j(\exp(-\nu_j)g_{j-1}A)$ is in the normal form.
\\
\\~

As a result we construct a sequence $\{B_j\}$ converging to $B\in N(G)$.
The matrix $ A $ is jet-by-jet equivalent to $ B $ and,
by Theorem \ref{Thm.Jet.by.Jet.equivalence.implies.Equivalence.For.$G^0$} it is $G^0$-equivalent to the normal form.
\\
\\
\\
2. To check uniqueness, let the matrices $B,\tB \in N(G) $ be $ G^0$-equivalent, i.e. $\tB=g.B$ for
$g=\exp\nu \in G^0$. Then $\pi_0\tB=\pi_0B$. Assume the equality $ \pi_{j-1}\tB = \pi_{j-1}B $ is proved.
Then $\pi_j(\tB-B)\in  Mat^{(j)}(m,n,p)$ and
\beq
\pi_{j-1}(g.B-B)=\pi_{j-1}(\tB-B) = 0.
\eeq
Hence, $ g \in \St_{j-1}(B)$, and $\nu\in\Ker\pi_{j-1}S_B'(0)$.
Further, by Lemma 3.1
\beq
\pi_j(\tB-B)=\pi_j(g.B - B)=\pi_jS_B'(0)\nu \in V^{(j)}_B.
\eeq
Since $ \tB, B $ are in the normal form, the inclusion $ \pi_j(\tB - B) \in
W^{(j)}(\pi_{j-1}B) $ holds. Hence $ \pi_j\tB = \pi_jB$ for any $j$. Hence, $\tB=B$.
\epr

Theorem \ref{Thm.Normal.Form.Exists} together with Theorem \ref{Thm.Finite.Determinacy.Equivalent.Condition}
 imply:
\bcor\label{Thm.Corollar.If.k.determined.Then.the.normal.form.is.Polyn}
If a matrix is $k$-determined then its normal form is a matrix of polynomials of degrees $\le k$.
\ecor
Indeed, the complementary spaces $W^{(j)}(A)$ will be zero for $j\ge k$.
\subsection{The inner product and the differential operators}
Let now $\K\subset\C$. Then we can apply the Euclidian structure to choose the complements $W^{(j)}(B)$.

Let $A(x)=\suml_{|I|\leq j}A_Ix^I$ and  $B(x) = \suml_{|I|\leq j}B_Ix^I$
be two polynomial matrices. We introduce the inner product
\beq
\bl A,B\br = \suml_{|I|\leq j}(A_I,B_I)I!, \hspace{1cm}I! = I_1! \ldots I_P!
\eeq
where $(A_I,B_I) = {\rm trace}B_I^*A_I$. In the space of the j'th jets of matrices, $Mat_j(m,m,p)\times Mat_j(n,n,p)$,
we introduce the similar inner product
\beq
\bl\nu,\mu\br=\bl\nu_l,\mu_l\br+\bl\nu_r,\mu_r\br
\eeq
for $\nu=(\nu_l,\nu_r)$, $\mu=(\mu_l,\mu_r)$.
 Then the projectors $ \pi_j $ are self-adjoint: $\bl \pi_j\nu,\mu\br = \bl\nu,\pi_j\mu\br$.
\\\\~

Given a formal matrix $B(x)=\suml_I B_Ix^I$, we introduce the formal differential operators
\beq
B^*\left(\frac{\di}{\di x}\right):=\suml_IB_I^*\frac{\di^{|I|}}{\di x^I}, \hspace{1cm}
(B^*)^T\left(\frac{\di}{\di x}\right):=\suml_I(B_I^*)^T\frac{\di^{|I|}}{\di x^I}.
\eeq
where $B^*_I:=\overline{B^T_I}$.

The first acts on the space of all polynomial $ m\times n $ matrices,
while the second on the space of all polynomial $ n\times m $ matrices.

The action of the differential operator $D_B$ on a polynomial matrix $P(x)$ is defined as
\beq\ber
D_BP:= \nu \in Mat(m,m,p)\times Mat(n,n,p)\\

\nu_l(x)=\left((B^*)^T\left(\frac{\di}{\di x}\right) P^T(x)\right)^T,\hspace{1cm}
\nu_r(x) = -B^*\left(\frac{\di}{\di x}\right) P(x).
\eer\eeq
\subsection{The normal form based on an inner product}\label{Sec.Normal.Form.from.Inner.Product} ~\\
Let $G$ be a Lie group and $\Lie(G)$ its Lie algebra. For each $j$ there is the natural inclusion:
$\Lie^{(j)}(G)\subset Mat^{(j)}(m,m,p)\times Mat^{(j)}(n,n,p)$.
Using the inner product we can define the orthogonal complement of $\Lie^{(j)}(G)$
and the orthogonal projection onto $\Lie^{(j)}(G)$.
Hence we have the collection of orthogonal projectors:
\beq
\gd_j: Mat^{(j)}(m,m,p)\times Mat^{(j)}(n,n,p)\rightarrow\Lie^{(j)}(G)
\eeq
\bex
For the Lie groups introduced above the projectors $\gd_j$ have a very simple form:
\beq\ber
\gd_j\nu = (\one-\pi_0)\nu \text{ for } G=G_{lr};\\
\gd_j\nu = (\one-\pi_0)(\nu_l,0) \text{ for } G=G_l;\\
\gd_j\nu = (\one-\pi_0)(0,\nu_r) \text{ for } G=G_r;\\
\gd_j\nu = (\one-\pi_0)\left(\frac{\nu_l+\nu_r}{2},\frac{\nu_l+\nu_r}{2}\right) \text{ for } G=G_c;\\
\gd_j\nu = (\one-\pi_0)\left(\frac{\nu_l-\nu_r^T}{2},\frac{\nu_r-\nu_l^T}{2}\right) \text{ for } G=G_T.
\eer\eeq
\eex
In the following we denote the projectors just by $\gd$, assuming that they act on the corresponding subspaces.
\bthe\label{Thm.Normal.Form.from.Inner.Product}
Every formal matrix $A(x)$ is $G^0$-equivalent to a unique matrix $B(x)$ whose homogeneous summands $B^{(j)}$
satisfy the equation
\beq\label{Eq.Normal.Form.With.Inner.Product}
\gd(D_BB^{(j)})(x)=\gd(D_B\pi_{j-1}f_j)(x), \ \ j = 0,1, \ldots
\eeq
with some formal matrices $f_j$.
\ethe
Hence, the subset $N(G)\subset Mat(m,n,p)$ of formal matrices $B$ satisfying the
condition (\ref{Eq.Normal.Form.With.Inner.Product}) is a normal form \wrt  $G$ and a canonical
one \wrt  $G^0$.
\\\bpr
Recall from Theorem \ref{Thm.Normal.Form.Exists} that for each $j$, having built $\pi_{j-1}(B)$,
 we should fix a complement to the vector space
$V^{(j)}(B)=\{\pi_j(\nu B)| \ \nu\in \Lie(G), \  \pi_{j-1}(\nu B)=0\}$.
 Recall also that $V^{(j)}(B)$ depends on the $(j-1)$'st jet of $B$ only, i.e.  $V^{(j)}(B)=V^{(j)}(\pi_{j-1}B)$.
Set
\beq
Mat^{(j)}(m,n,p)\supset W^{(j)}(B):=(V^{(j)}(B))^\perp=\Big(\pi_jS_B'(0)(\Ker \pi_{j-1}S_B'(0))\Big)^\perp
\eeq
where $\perp$ means the orthogonal complement.
It suffices to show that the matrix $B$ lies in the set $N(G)$, i.e. $\forall j\in\mN:$ $B^{(j)}\in W^{(j)}(B)$,
if and only if equation (\ref{Eq.Normal.Form.With.Inner.Product}) holds.

Consider the linear map
\beq\ber
Mat(m,m,p)\times Mat(n,n,p)\stackrel{\B}{\rightarrow}Mat(m,n,p)
\\
\hspace{1.5cm}\nu=(\nu_l,\nu_r)\mapsto \B \nu= \nu_lB - B\nu_r
\eer\eeq
Then, for $\nu\in\Lie(G)$ have $\pi_jS_B'(0)\nu = \pi_j\B \gd\nu$.
\\
\\\\
We claim that $B^{(j)}\in W^{(j)}(B)$ iff
\beq\label{Eq.Inside.Proof.1}
(\pi_j\B\gd)^*B^{(j)} \in Im(\pi_{j-1}\B \gd)^*
\eeq
Indeed, suppose (\ref{Eq.Inside.Proof.1}) holds. Then $(\pi_j\B \gd)^*B^{(j)} \in (\Ker \pi_{j-1}\B \gd)^\perp$,
implying
\beq\label{Eq.Inside.Proof.2}
B^{(j)} \in (\pi_j\B \gd)(\Ker \pi_{j-1}\B \gd)^\perp =W^{(j)}(B)
\eeq
Conversely, suppose (\ref{Eq.Inside.Proof.2}) is valid. Then
\beq
<B^{(j)},\pi_j\B \gd\nu> = 0, \nu \in \Ker \pi_{j-1}\B \gd.
\eeq
Therefore
\beq
(\pi_j\B \gd)^*B^{(j)} \in(\Ker \pi_{j-1}\B \gd)^\perp ={\rm Im}(\pi_{j-1}\B \gd)^*,
\eeq
i.e. (\ref{Eq.Inside.Proof.1}) holds.
\\
\\
\\
In order to finish the proof it remains to compute the conjugate maps $ (\pi_j\B \gd)^* $ and $ (\pi_{j-1}\B \gd)^*$.
Note that $\pi_j\B \gd\nu =\pi_j\gd(\nu_l) B- \pi_j B\gd(\nu_r)$ for $\nu=(\nu_l,\nu_r)$.
We claim that
\beq
(\pi_j\B \gd)^*P = \gd\nu=\Big(\gd(B^{*T}P^T)^T,\gd B^*P\Big)
\eeq
To check this we compute $\bl Q,(\pi_j\B \gd)^*P\br$ for $Q=\sum Q_I x^I$.
\beq\ber
\bl Q,(\pi_j\B \gd)^*P\br=\bl \pi_j\B \gd Q,P\br=\bl \sum (Q_l)_J B_Ix^{I+J}-\sum B_I(Q_r)_Jx^{I+J},\pi_j(P)\br=
\\=\bl Q_l,B^*(\frac{\partial}{\partial x}) \pi_j P \br-\bl Q_r,\Big(B^*(\frac{\partial}{\partial x})^T \pi_j P^T \Big)^T\br
=\bl Q,(\B\gd)^*\pi_j P\br
\eer\eeq

As a result,
\beq
(\pi_j\B \gd)^*P = \gd D_B\pi_j, \ (\pi_{j-1}\B \gd)^*P = \gd D_B\pi_{j-1},
\eeq
proving the statement.\epr
\subsection{Corollaries and examples.}\label{Sec.Corollaries.And.Examples}
Let $\Ak(x) $ be a homogeneous matrix of degree $k$ and let
$A(x)=\Ak(x) + {\rm terms\ of\ orders}\ \geq k+1$.

Then the biggest powers in equation (\ref{Eq.Normal.Form.With.Inner.Product}) for the above listed groups
are equal to $k-j$, and arise only in the left side. If $G=G_c$ then the same is true for the matrix
\beq
A(x) = \gl \one +\Ak(x) \ldots,  \ \gl \in \K.
\eeq
Taking into account the structure of the projectors $ \gd, $ we arrive at
\\
\bcor Let $A(x)=\Ak(x)+(\text{terms of orders}\geq k+1)$. Then
\\1. The matrix $ A $ is left equivalent to a matrix $\Ak(x) + b(x) $ satisfying
\beq\label{Eq.Corollary.Eq1}
(\Ak^*)^T\left(\frac{\di}{\di x}\right) b^T(x) = 0
\eeq
2. The matrix $ A $ is right equivalent to a matrix $\Ak(x)+b(x)$ satisfying
\beq\label{Eq.Corollary.Eq2}
\Ak^*\left(\frac{\di}{\di x}\right) b(x)=0
\eeq
3. The matrix $A$ is two-sided equivalent to a matrix $\Ak(x)+b(x)$ satisfying both of
the relations (\ref{Eq.Corollary.Eq1}) and (\ref{Eq.Corollary.Eq2}).
\\4. If $m=n$, then the matrix $ A $ is congruent to a matrix $\Ak(x)+b(x)$ satisfying
\beq
(\Ak^*)^T\left(\frac{\di}{\di x}\right) b^T(x) + \Ak^*\left(\frac {\di}{\di x}\right) b(x) = 0.
\eeq
5. If $ m=n, $ then every matrix $ \gl \one + A,\  \gl \in \K $ is conjugate to a matrix $ \gl \one + P + b $ satisfying
\beq
\left((\Ak^*)^T\left(\frac{\di}{\di x}\right) b^T(x)\right)^T =\Ak^*\left(\frac{\di}{\di x}\right) b(x).
\eeq
\ecor
The statement 3 proves Theorem \ref{Thm.Normal.Form.Main.Theorem} from the Introduction.

The relation similar to Corollary \ref{Thm.Introduction.Corollary.Minimal.Resolution} for the group
$G_{lr}$ from the Introduction takes on the form
\beq\ber
\gL^*b(x) = 0 \text{ for the group } G_l\\
b(x)\gL^* = 0 \text{ for the group } G_r\\
\gL^*b(x) = b(x)\gL^* \text{ for the conjugacy }\\
(\gL^*)^Tb^T(x) + b(x)\gL^* = 0 \text{ for the  congruence}
\eer\eeq

In addition to Corollary 1.2 from Introduction we obtain
\bcor
Let $ m=n $ and $\Ak(x) = {\rm diag}\big(l_1(x) \ldots l_m(x)\big)$ with homogeneous polynomials of degree $ k. $ Then
\\1. Every matrix $ A(x) =\Ak(x) +(\text{terms of orders}\geq k+1)$ is left equivalent to a
matrix $\Ak(x)+b(x)$ satisfying $l_j^*\left(\frac{\di}{\di x}\right) b_{ij}(x) = 0$.
\\2. The matrix $A(x)$ is right equivalent to a matrix $\Ak(x)+b(x)$ satisfying
$l_i^*\left(\frac{\di}{\di x} \right) b_{ij}(x) = 0$.
\\3. If $ m=n, $ then the matrix $ A(x) $ is congruent to a matrix $\Ak(x)+b(x)$ satisfying
$l_I^*\left(\frac{\di}{\di x}\right) b_{ij}(x) +l_j^*\left(\frac{\di}{\di x}\ \right) b_{ij}(x) = 0$.
\\4. If $ m=n, $ then every matrix $ \gl \one + A(x),\  \gl \in \K $ is conjugate to a
 matrix $ \gl \one +\Ak(x) + (b_{ij}(x)) $ satisfying
$\left( l_i^*\left(\frac{\di}{\di x}\right) - l_j^*\left(\frac{\di}{\di x}\right)\right)b_{ij}(x) = 0$.
\ecor
\bex Let $m=n$ and suppose the polynomials $l_i$ are linear, i.e. $k=1$ and
$l_i(x) = \suml_{s=1}^p\ga_{is}x_s$, for $i = 1, \ldots, m$.
Then Corollary 1.2 for the group $ G_{lr} $ gives
\beq
\suml_{s=1}^p\overline{\ga}_{is}\frac{\di b_{ij}(x)}{\di x_s} = 0, \hspace{1cm}
 \suml_{s=1}^p\overline{\ga}_{is}\frac{\di b_{ij}(x)}{\di x_s} = 0.
\eeq
In particular, if $ p=2 $ and the functionals $ l_s $ are pair wise non-colinear,
then $ b_{ij} = 0 $ for $ i \neq j. $ Hence, the normal form is a diagonal matrix.
\eex
\bex Let $m=n$ and $\Ak(x)=l_1(x)\one+l_2(x)J$. Here $l_1,l_2$ are homogeneous polynomials of degree $k$, while
\beq
\tiny J=\bpm 0&1&.&0\\.&.&.&.\\0&.&0&1\\0&.&.&0\epm
\eeq
is the nilpotent Jordan block. Then the relations of Corollary \ref{Thm.Introduction.Corollary.For.Diag.Matrix}
take on the form
\beq
l_1^*\left(\frac{\di}{\di x}\right) b_{1j}(x) = 0, \ \
l_1^*\left(\frac{\di}{\di x}\right) b_{ij}(x) +l_2^*\left(\frac{\di}{\di x}\right) b_{i-1j} = 0
\eeq
for $j=1,2, \ldots, m$ and $i\geq 2$. Besides,
\beq
l_1^*\left(\frac{\di}{\di x}\right) b_{im}(x) = 0, \ \
l_1^*\left(\frac{\di}{\di x}\right) b_{ij}(x)+l_2^*\left(\frac{\di}{\di x}\right) b_{ij+1}=0
\eeq
for $ i = 1,2, \ldots, m $ and $ j \leq m-1$. It follows that
\beq\label{Eq.Last.Corollary}
l_1^*\left(\frac{\di}{\di x}\right) b_{ij}(x) = 0, \hspace{1cm}
l_2^*\left(\frac{\di}{\di x}\right) b_{ij}(x)=0
\eeq
for $ j \geq i+1. $ Hence, every matrix
\beq
A(x) = l_1(x)\one + l_2(x)J + {\rm terms\ of\ orders\ } \geq k+1
\eeq
is two-sided equivalent to a matrix $ B(x) = l_1(x)\one + l_2(x)J +
(b_{ij}(x)) $ where $ b_{ij} $ satisfy (\ref{Eq.Last.Corollary}).
\eex
\bex
In particular, consider the case of two variables, $p=2$, and assume $l_i(x)$ are linear and linearly independent.
Then equation (\ref{Eq.Last.Corollary}) implies $b_{ij} = 0$ for $j \geq i+1$.
Hence, in this case the matrix $ A $ is two-sided equivalent to a matrix of the form
\beq
l_1(x)\one + l_2(x)J +
\bpm ~* & 0 & 0 & \dots & 0 \\ * & * & 0 &\dots& 0 \\.&.&.&\dots&.\\ * & * & * &\dots& * \epm
\eeq
By the linear change of coordinates we can make $l_1(x)=x_1$ and $l_2(x)=x_2$, then the entries of
the low triangular matrix satisfy:
\beq
\di_x b_{ij}+\di_y b_{i-1,j}=0,\hspace{1cm} \di_x b_{ij}+\di_y b_{i,j+1}=0
\eeq
This implies $\di_y b_{i-1,j}=\di_y b_{i,j+1}$. Combining with the last equations it gives
$\di_x b_{i+1,j}=\di_x b_{i,j-1}$. As $b_{ij}$ is of order at least two we get: $b_{i,j}=b_{i+1,j+1}=:\gga_{i-j}$.
Finally we obtain that $A$ is two-sided equivalent to
\beq
x\one+yJ+\bpm ~\gga_0 & 0 & 0 & \dots & 0 \\ \gga_1 & \gga_0 & 0 &\dots& 0 \\.&.&.&\dots&.
\\ \gga_{n-1} & \gga_{n-2} & . &\dots& \gga_0 \epm,\hspace{1cm} \di_x \gga_i+\di_y\gga_{i-1}=0
\eeq
\eex
\appendix
\section{Dependence on the choice of the base ring}\label{Sec.Dependence.On.The.Base.Ring}
Let the field $\K$ have a non-trivial valuation, so that the convergence of a
series is defined. One can consider matrices whose entries are
\li formal series- $\K[[x_1,..,x_p]]$ or
\li locally converging series- $\K\{x_1,..,x_p\}$ or
\li rational functions that are regular at the origin- $\K[x_1,..,x_p]_{(\frm)}$,
(i.e. fractions of polynomials, whose denominators do not
vanish at the origin). Here $\frm\subset\K[x_1,..,x_p]$ is the maximal ideal.

Correspondingly we have the notions of formal/locally converging/rational-G-equivalences,
 formal/locally converging/rational-G-determinacy etc.
Two questions occur naturally:

(injectivity) Let $A_1$, $A_2$ be matrices with rational/locally converging entries. Suppose they are
formally-G-equivalent. Are they rationally/locally converging-G equivalent?

(surjectivity) Which formal matrices are G-equivalent to matrices with locally converging/rational entries?

In the discussion below many things are well known in Commutative Algebra (e.g. \cite{Eisenbud-book},
\cite{Yoshino-book}), but they seem to be less known in other areas.

Recall that even for a formal series $f\in\K[[x_1,..,x_p]]$ one can speak about the
corresponding hypersurface $\{f=0\}\subset(\K^p,0)$, its singularities, local irreducibility etc.
 Though the series might not converge off the origin.
\subsection{Injectivity}
We consider here only those subgroups $G\subset G(m,n,p)$ that are defined by polynomial equations in matrix entries.
More precisely $(U,V)\in G$ iff the entries of the matrices $(U,V)$ satisfy a finite collection of polynomial
equations with constant coefficients. (For example $G_l$,$G_r$,$G_{lr}$,$G_c$,$G_T$ are such subgroups.)
\bthe
1. Let $A_1$,$A_2$ be matrices with locally converging entries, let $G\subset G(m,n,p)$. If $A_1$,$A_2$
are formally-G-equivalent, then they are locally convergent-G-equivalent.
\\2. Let $A_1$, $A_2$ be matrices with rational entries. Let $G$ be a subgroup of $G(m,n,p)$ defined
by equations linear in matrix entries. For example $G$ is one of $G_l$,$G_r$,$G_{lr}$,$G_c$.
If $A_1$,$A_2$ are formally-G-equivalent, then they are rationally-G-equivalent.
\ethe
\bpr
1. If $A_1$,$A_2$ are formally-G-equivalent then $UA_1=A_2V$, where $U$, $V$
are invertible at the origin, and satisfy some additional algebraic conditions
(depending on $G$).
So, if the entries of $A_1,A_2$ are locally converging, then by Artin approximation theorem \cite{Artin68}
the matrices $U,V$ can be chosen with locally converging entries.
\\
\\
2. Suppose $A_1,A_2$ have rational entries and $G$ satisfies the assumption.
Then the conditions on $(U,V)$ are:
\li {\em linear} equations for the entries of $U,V$. These arise from $UA_1=A_2V$ and from the defining
conditions of the group (e.g. for $G_l$: $V=\one$, for $G_c$: $U=V$).
The coefficients in these equations are constants or the entries of $A_1,A_2$, i.e. rational functions.
\li the {\em non-degeneracy} condition:  $U$ and $V$ are invertible at the origin.

Note that some of the linear equations above can be non-homogeneous, e.g. $V=\one$ for $G_l$. So the set of
all the pairs of matrices satisfying these linear conditions is an affine space, but in general not a linear one.
Hence we 'homogenize' the equations by introducing new variables.
For example for $G_l$ replace $V=\one$ by $V=\tV$, with $\tV$ the matrix whose entries are additional variables,
 at the end we will impose the additional constraint $\tV=\one$. Note that the equations of these
 additional constraints are linear, with constant coefficients.

Let
$E\subset Mat\big(m\times m,\K[x_1,..,x_p]_{(\frm)}\big)\oplus Mat\big(n\times n,\K[x_1,..,x_p]_{(\frm)}\big)\oplus..$
 be the set of all the tuples $(U,V,additional\ variables)$, whose entries are rational functions regular at the origin,
 such that the tuple satisfy the homogenized linear conditions as above. So $E$ is
a vector space. In fact $E$ is a module over $\K[x_1,..,x_p]_{\frm}$ with the action $f(U,V,..):=(fU,fV,..)$.

Let $E_{formal}\subset Mat\big(m\times m,\K[[x_1,..,x_p]]\big)\oplus Mat\big(n\times n,\K[[x_1,..,x_p]]\big)\oplus..$
be the set of all the tuples $(U,V,..)$, with formal entries, satisfying the homogenized linear
conditions as above. So $E_{formal}$ is a module over $\K[[x_1,..,x_p]]$.
We claim that $E_{formal}=\K[[x_1,..,x_p]] E$, i.e. one can choose a basis of $E_{formal}$ consisting of
the elements of $E$.

First, observe that $E$ is a finitely generated module over $\K[x_1,..,x_p]_{(\frm)}$.
For example, consider the ideal in $\K[x_1,..,x_p]_{(\frm)}$
generated by $U_{1,1}$ for all $(U,V,..)\in E$.
By Hilbert basis theorem \cite[pg.371]{Eisenbud-book} this ideal has a finite basis, say $\{(U_i,V_i,..)\}_i$.
Correspondingly, the module $E$ decomposes: $E=E'\oplus Span\big(\cup_i(U_i,V_i,..)\big)$. Here
$Span\big(\cup_i(U_i,V_i,..)\big)$ is the submodule of $E$ generated by $\{(U_i,V_i,..)\}_i$, while $E'$
is the submodule generated by those $(U,V,..)$ that have $U_{11}=0$. Continue in this way over all the
entries of $U$ and $V$, to get a finite basis for $E$.

Next, consider the free injective resolution of $E$:
\beq
0\to E\to Mat\big(m\times m,\K[x_1,..,x_p]_{(\frm)}\big)\oplus Mat\big(n\times n,\K[x_1,..,x_p]_{(\frm)}\big)\oplus..
\stackrel{\phi}{\to}...
\eeq
Here $Mat\big(m\times m,\K[x_1,..,x_p]_{(\frm)}\big)\oplus Mat\big(n\times n,\K[x_1,..,x_p]_{(\frm)}\big)\oplus..$
is considered as a free module over $\K[x_1,..,x_p]_{(\frm)}$.
The map $\phi$ corresponds to all the linear homogenized equations imposed on $(U,V,..)$.

Now take the completion $\K[x_1,..,x_p]_{(\frm)}\to \K[[x_1,..,x_p]]$.
 As the completion functor is exact, \cite[pg. 198]{Eisenbud-book}, the resolution is preserved:
\beq
0\to \K[[x_1,..,x_p]]E\to Mat\big(m\times m,\K[[x_1,..,x_p]]\big)\oplus Mat\big(n\times n,\K[[x_1,..,x_p]]\big)\oplus..
\stackrel{\phi}{\to}...
\eeq
But the last row is the resolution of $E_{formal}$. Hence $E_{formal}=\K[[x_1,..,x_p]]E$.
\\
\\
Now impose the additional conditions on the new variables introduced to homogenize the initial conditions
 (e.g. $\tV=\one$, for $G_l$). As they are all linear, with constant coefficients,
 we still have the property: if $(U,V)$ is a formal solution of the initial linear equations,
  then $(U,V)=\sum f_i (U_i,V_i)$ for some rational solutions $(U_i,V_i)$ and $f_i\in \K[[x_1,..,x_p]]$.
\\
\\
Finally, suppose $UA_1=A_2V$ has a formal solution for $(U,V)\in G$, in particular $(U,V)$ are invertible
at the origin.
 Hence, $U':=\sum jet_0(f_i)U_i$ and $V':=\sum jet_0(f_i)V_i$ are locally invertible matrices
of rational functions satisfying $U'A_1=A_2V'$.
\epr
\beR
The second statement of the proposition is not true for $G=G_T$ as in this case the conditions on $U,V$ are non-linear.
For example, let $A_2=(1+x)A_1$ be $1\times 1$ matrices, i.e. functions.
Then for $A_1=UA_2U^T$ one has $U^2=1+x$, i.e. $U$ cannot be rational.
\eeR
\bcor
Let the matrix $A$ have locally converging entries. Suppose $A$ is formally-finitely-$G$-determined, i.e.
the conditions of Theorem \ref{Thm.Finite.Determinacy.Equivalent.Condition} or Corollary
\ref{Thm.Corollary.Finite.Determinacy} are satisfied.  Then, $A$ is locally-converging-finitely-$G$-determined
and $A$ is locally-converging-$G$-equivalent to a matrix of polynomials.
Further, by  Corollary \ref{Thm.Corollar.If.k.determined.Then.the.normal.form.is.Polyn}
the normal form of $A$ is polynomial.
\ecor
\subsection{Surjectivity}
Most matrices with formal/locally convergent entries are not $G_{lr}$-equivalent to locally convergent/rational matrices.
\bex 
Let $A(x,y) =y-xf(x)$ be a ``$ 1\times 1 $ matrix'' of two variables, where $f(x)$ is a
formal (but not locally converging) series or a locally converging series (but not a rational function).

Assume there exist a formal $1\times 1$ matrix $U(x,y)$, invertible at the origin, such that $U(x,y)A(x,y)$
is a locally converging series/a rational function. Note that if $U(x,y)A(x,y)$ vanishes at some point then
$A(x,y)$ vanishes too.

Hence if $U(x,y)A(x,y)$ is a rational function then $\{U(x,y)A(x,y)=0\}\subset\K^2$ is
an algebraic curve, which is defined also as $\{y-xf(x)=0\}$. Implying that $f(x)$ is rational, contradiction.

Similarly, if $U(x,y)A(x,y)$ is a locally convergent power series, then it defines a locally analytic curve.
On this curve $y=xf(x)$, i.e. $f(x)$ must be convergent at every point of this curve, contradiction.
\eex
An immediate necessary condition for a square matrix to be equivalent to a matrix of locally convergent series/rational
functions is:
 $\det(A)$ is a locally convergent series/rational function, up to an invertible factor.
 Or, the ideal $\bl\det(A)\br\subset\K[[x_1,..,x_p]]$ is generated by a locally convergent series/rational function.

A stronger condition: let $I_k(A)$ be the ideal in $\K[[x_1,..,x_p]]$ generated by all the $k\times k$
minors of $A$. Note that these ideals are invariant under $G_{lr}$ equivalence.
Hence, if $A$ is equivalent to a matrix of locally convergent series/rational functions,
then all the ideals $I_k(A)\subset\K[[x_1,..,x_p]]$ are generated by locally convergent series/rational functions.
All of these conditions are relevant, as the following example shows.
\bex
Consider the matrix with entries in $\K[[x,y,z,q,w]]$:
\beq
A=\begin{pmatrix} z&y+x^2f_1(x)&0\\0&w&x+y^2f_2(y)\\0&0&q\end{pmatrix}
\eeq
It has a polynomial determinant and the ideal of its entries is the maximal ideal,
  $I_1(A)=<x,y,z,w,q>$. In particular this ideal is polynomially generated.
If $f_1(x)$, $f_2(x)$ are formal but not locally convergent/locally convergent but not rational, then
 $A$ is not $G_{lr}$ equivalent to a matrix with locally convergent/rational entries.
Because $I_2(A)$ is not generated by locally convergent/rational elements.
\eex
\beR
As has been proved recently, \cite[proposition 1.6]{Keller-Murfet-Van den Bergh2008},
for any formal matrix $A\in Mat(m,m,p)$, with arbitrary field $\K$, there exists a matrix $B\in Mat(n,n,p)$
such that $A\oplus B$ is $G_{lr}$ equivalent to a matrix of rational functions.
\eeR
\subsection{The case of two variables}
Every formal matrix of one variable is $G_{lr}$ equivalent to a polynomial matrix, e.g. see the normal form
\S\ref{Sec.One.Variable.Case}. We prove that to some extent this is true in the case of two variables.

In this section $\K$ is algebraically closed, $A$ is an $m\times m$ matrix, with entries in $\K[[x,y]]$.
We always assume $A|_0=\zero$ (cf. Corollary \ref{Thm.Introduction.Corollary.Minimal.Resolution})
 and $\det(A)\not\equiv0$.
In addition we assume: the plane curve $C:=\{det(A)=0\}$ is {\em reduced}, i.e. it has no multiple components.
(Note though that $C$ can be reducible.)

First, recall the situation with functions.
\bprop\label{Thm.Formal.Function.Is.Polynomial.Mod.A.Curve}
Let $f,g\in\K[[x,y]]$ be relatively prime. Then $g=u g'\mod(f)$, where $u\in\K[[x,y]]$ is invertible and $g'\in\K[x,y]$.
\eprop
Mote generally, let $h=(f,g)$ be the greatest common divisor, then $g=u h g'\mod(f)$ with $u$
invertible and $g'$ a polynomial.
This property is well known, but we could not find a reference. Hence we give a proof.
\\\bpr
Step 1. First suppose the formal curve $C=\{f=0\}\subset(\K^2,0)$ is locally irreducible. Let $\tC\stackrel{\nu}{\to}C$ be the
normalization, i.e. a smooth curve germ and a finite morphism that is an isomorphism outside the singularity of $C$.
(See e.g. \cite[pg.125-129]{Eisenbud-book}.)
This corresponds to the embedding of the local rings: $\quotient{\K[[x,y]]}{(f)}\stackrel{\nu^*}{\hookrightarrow}\K[[t]]$.
By the finiteness of the morphism, the quotient $\quotient{\K[[t]]}{\K[[x,y]]/(f)}$ is a finite dimensional vector space.

The normalization induces the valuation
\beq
val:\K[[x,y]]\to\mN,\ \ g\to val(g):=ord_t\nu^*(g)= ord_t(g(x(t),y(t)))
\eeq
By the finiteness of the quotient above, there exists the conductor, i.e.
the minimal number $c\in\mN$ such that any bigger number $d>c$
is realized as the valuation of some function: $d=val(g)$, for $g=\K[[x,y]]$.
\\
\\
Step 2. Let $g\in \K[[x,y]]$, not a polynomial. Then can decompose $g$ into the sum of a polynomial
and some series of high valuation: $g=g_{pol}+g_{high}$, where $val(g_{high})>c+val(g_{pol})$.
By the existence of conductor, there exists $h\in\K[[x,y]]$ such that $val(h)=val(g_{high})-val(g_{pol})$,
i.e. $val(h g_{pol})=val(g_{high})$. Then for some number $\gb$ one has: $val(hg_{pol}+\gb g_{high})>val(g_{high})$.
Hence $(1+h)g=g_{pol}+g'_{high}$ with $val(g'_{high})>val(g_{high})$.

Continue this process, to get in the limit: $g=u g_{pol}+g_{high}$ such that $val(g_{high})=\infty$ and $u$
is invertible. But then $g_{high}$ is divisible by $f$. Hence the statement.
\\
\\
Step 2'. If the curve $C=\{f=0\}$ is locally reducible, $C=\cup^r_{i=1} C_i$, then the normalization $\coprod \tC_i\to\cup C_i$
induces the multi-valuation $\{val_i\}:\K[[x,y]]\to\mN^{\oplus r}$. The quotient of the local ring is
still a finite dimensional vector space, hence the conductor still exists. Continue as above, using
that $g,f$ are mutually prime.
\epr
\bthe
Let $A\in Mat(m,m,2)$ be a formal matrix in two variables. Assume $\det(A)$ is locally
convergent (up to an invertible factor).
Then $A$ is $G_{lr}$ equivalent to a matrix with locally converging entries.
\ethe
\bpr
{\em Step 1.} We can consider $A$ as a matrix with entries in $\quotient{\K[[x,y]]}{(\det A)}$.
(Recall that $\det A\not\equiv0$ and is reduced, and $A|_0=\zero$.)
Namely we consider $A\mod(\det A)$. The $G_{lr}$ equivalence over $\K[[x,y]]$ descends to that over
$\quotient{\K[[x,y]]}{(\det A)}$.

Conversely, if $A\sim B$ over $\quotient{\K[[x,y]]}{(\det A)}$ then they are
equivalent over $\K[[x,y]]$. Indeed, suppose $A=UBV+(\det A)Q$, where $U,V$ are invertible and $Q$ is some formal matrix.
Recall that $\det(A)=A\cAv$, where $\cAv$ is the adjoint matrix. Thus we get $A(\one-\cAv Q)=UBV$. As $A|_0=\zero$
we get that $(\cAv Q)|_0=\zero$, hence the matrix $(\one-\cAv Q)$ is invertible. Therefore:
$A=(\one-\cAv Q)^{-1}UBV$.

So, it is enough to show that $A$  is equivalent to a locally convergent matrix modulo $\det A$.
In fact, we will show this for $\cAv$ and then achieve the statement for $A$ too. Note that
an equivalence transformation $A\to UAV$ results in the equivalence $\cAv\to V^{-1}\cAv U^{-1}$.
\\
\\
{\em Step 2.} From now on consider $A$ and $\cAv$ modulo $\det A$. By the previous proposition
we can assume $\cAv$ in the form $\{u_{ij}g_{ij}\}$, where $u_{ij}$ are invertible and $g_{ij}$
are locally converging. Consider the equivalence transformation
\beq
\cAv\to U\cAv V, \ \ U=\bpm u^{-1}_{11}&0&..&0\\0&u^{-1}_{21}&..&0\\..&..&\\0&..&0&u^{-1}_{m1}\epm,\ \
V=\bpm u^{-1}_{11}&0&..&0\\0&u^{-1}_{12}&..&0\\..&..&\\0&..&0&u^{-1}_{1m}\epm
\eeq
So the first row and column of $U\cAv V$ have locally convergent entries.  From now on we assume $\cAv$ in this form.
\\
\\
{\em Step 3.} By definition of $\det$, the matrix $A$ is degenerate, when restricted to the curve $C=\{\det A=0\}$.
 At a point $pt\in C$ the corank of $A$ is not bigger than the multiplicity $mult(C,pt)$.  Hence, as $C$
 is reduced, the corank of $A$ at the smooth points of $C$ is one.

Recall that $A\cAv|_C=\zero|_C$, hence at the smooth points of $C$ the rank of $\cAv$ is one. Namely,
any two rows/columns are dependent.  Thus in particular, for any entry $\cAv_{ij}$ one has:
$\cAv_{ij}=\frac{\cAv_{i1}\cAv_{1j}}{\cAv_{11}}\in\quotient{\K[[x,y]]}{(\det A)}$. Note that the right hand side
 is {\em locally convergent}.

Hence, the above equivalence transformation results in a locally convergent matrix $\cAv$.
Thus $A\sim(\cAv)^\vee$ is locally convergent too.
\epr

\bthe
Suppose $A$ is a formal matrix in two variables and its determinant is a non-zero, reduced
polynomial, up to an invertible factor.
\\ * If $\det A$ is irreducible (as a formal series) then $A$ is $G_{lr}$-equivalent to a
polynomial matrix.
\\ * More generally, suppose the decomposition of $\det A$ into irreducibles is the same over
formal series and over rational functions. Namely, in the decomposition into irreducible formal factors
$\det A=f_1\cdot\cdot\cdot f_k$ all $\{f_i\}$ can be chosen as polynomials (up to
multiplication by an invertible). Them $A$ is $G_{lr}$ equivalent to a polynomial matrix.
\ethe
\bpr
The proof goes precisely as in the locally convergent case. After {\em Step 2.} we have a matrix $\cAv$ whose
first row and column are polynomials. Hence in {\em Step 3.} we get
$\cAv_{ij}=\frac{\cAv_{i1}\cAv_{1j}}{\cAv_{11}}$, where on the right we have a rational function and on the
left a regular function. So $\cAv_{ij}$ is a regular rational function, i.e. a fraction of two polynomials,
with non-vanishing denominator.

Finally, multiply $\cAv$ by all such (invertible) denominators, to get: $\cAv$ is a polynomial matrix. From here obtain
that $A$ is polynomial too.
\epr

\end{document}